\documentclass[a4paper]{article}
\usepackage[latin1]{inputenc}
\usepackage[T1]{fontenc}
\usepackage[francais,english]{babel}
\usepackage{amssymb,color}
\usepackage{amsmath,amsfonts,amsthm}
\usepackage{amsfonts}
\usepackage{amssymb}
\usepackage{latexsym}
\usepackage{comment}
\usepackage{ifthen}
\usepackage{hyperref}

\def\R{\mathbb R}
\def\C{\mathbb C}

\newtheorem{thm}{Theorem}
\newtheorem{Prop}{Proposition}

\newtheorem{rk}{Remark}

\newboolean{highlight}
\setboolean{highlight}{true}
\definecolor{purple}{rgb}{0.89,0.15,0.99}

\author{
Karine \textsc{Beauchard}
\footnote{CMLS, Ecole Polytechnique, 91128 Palaiseau Cedex, France,
email: Karine.Beauchard@math.polytechnique.fr (corresponding author)},
Jean-Michel \textsc{Coron}
\footnote{Sorbonne Universit\'{e}s, UPMC Univ Paris 06, UMR 7598, Laboratoire
Jacques-Louis Lions, F-75005, Paris, France
email: coron@ann.jussieu.fr},
Holger \textsc{Teismann}\footnote{
Department of Mathematics and Statistics, Acadia University, Wolfville, NS, Canada
email: hteisman@acadiau.ca}
}

\title{Minimal time for the bilinear control of Schr\"{o}dinger equations}
\date{}

\begin{document}

\maketitle

\begin{abstract}
We consider a quantum particle in a potential $V(x)$ ($x \in \mathbb{R}^N$) 
subject to a (spatially homogeneous) time-dependent electric field $E(t)$, which plays the role of the control.
Under generic assumptions on $V$, this system is approximately controllable on the $L^2(\mathbb{R}^N,\mathbb{C})$-sphere,
in sufficiently large times $T$, as proved by Boscain, Caponigro, Chambrion and Sigalotti \cite{Boscain_CCS}.
In the present article, we show
that this approximate controllability result is false in small time.
As a consequence, the result by Boscain et al. is, in some sense, optimal 
with respect to the control time $T$.
\end{abstract}

\section{Introduction}

\subsection{Main result}

In this article, we consider quantum systems whose dynamics can be described by a linear Schr\"{o}dinger equation of the form
\begin{equation} \label{Schro}
\left\lbrace \begin{array}{ll}
\displaystyle i \partial_t \psi(t,x) = \Big( - \frac{1}{2} \Delta  + V(x) - \langle E(t) , x \rangle \Big) \psi(t,x)\,, & (t,x) \in (0,T)\times\mathbb{R}^N\,, \\
\displaystyle \psi(0,x)=\psi_0(x)\,,                                                                                               & x \in \mathbb{R}^N\,.
\end{array}\right.
\end{equation}
Here, $N \in \mathbb{N}^*$ is the space dimension,
$\langle.,.\rangle$ is the usual scalar product on $\mathbb{R}^N$, 
$V:x \in \mathbb{R}^N  \rightarrow  \mathbb{R}$,
$E:t \in (0,T) \rightarrow \mathbb{R}^N$ and
$\psi:(t,x) \in (0,T)\times\mathbb{R}^N \rightarrow \mathbb{C}$
are a static potential, a time-dependent  electric field, and the wave function, respectively.
This equation represents a quantum particle in the potential $V$ subject to the electric field $E(t)$.  
Planck's constant and the particle mass have been set to one.

System (\ref{Schro}) is a control system in which
the state is the wave function $\psi$, that belongs to the unitary $L^2(\mathbb{R}^N,\mathbb{C})$-sphere, denoted by $\mathcal{S}$;
and the control is the electric field $E$. 
Such systems have applications in modern technologies such as Nuclear Magnetic Resonance, quantum chemistry and quantum information science.
The expression 'bilinear control' refers to the `bilinear' nature w.r.t. $(E,\psi)$ of the term $\langle E(t) , x \rangle \psi$.


We are interested in the minimal time required to achieve approximate controllability of system \eqref{Schro}.
Since in (\ref{Schro}) decoherence is neglected, in realistic scenarios  the model may only
be applicable for small times $t$ (typically on the order of several periods of the ground state).
Thus, to be practically relevant,
controllability results need to be valid
for time intervals for which equation (\ref{Schro}) remains a reasonable model.
Therefore quantification of the minimal control time is an important issue.
\\

First, we recall a classical well-posedness result \cite{Fujiwara}, which we quote from \cite{carles2011}. 
We consider potentials $V$ that are smooth and subquadratic, i.e.
\begin{equation} \label{Ass_V}
V \in C^\infty(\mathbb{R}^N)\, \text{ and, }\,
\forall \alpha \in \mathbb{N}^N  \text{ such that }  |\alpha|\geqslant 2,\,  \,\partial_x^\alpha V \in L^\infty(\mathbb{R}^N)\,.
\end{equation}

\begin{Prop} \label{Prop:WP}
Consider $V$ satisfying assumption (\ref{Ass_V}) and $E \in L^\infty_{loc}(\mathbb{R},\mathbb{R}^N)$.
There exists a strongly continuous map $(t,s) \in \mathbb{R}^2 \mapsto U(t,s)$,
with values in the set of unitary operators on $L^2(\mathbb{R}^N,\mathbb{C})$, such that
$$U(t,t)=\text{Id}\,, \quad U(t,\tau) U(\tau,s)=U(t,s)\,, \quad U(t,s)^*=U(s,t)^{-1}\,, \quad \forall t,\tau,s \in \mathbb{R}$$
and for every $t,s \in \mathbb{R}$, $\varphi \in L^2(\mathbb{R}^N,\mathbb{C})$, the function $\psi(t,x):=U(t,s)\varphi(x)$ solves
the first equation of (\ref{Schro}) with initial condition $\psi(s,x)=\varphi(x)$.
\end{Prop}

For $V$ satisfying (\ref{Ass_V}), we introduce the operator
$$D(A_V):=\left\{ \varphi \in L^2(\mathbb{R}^N) ; -\Delta \varphi + V(x) \varphi \in L^2(\mathbb{R}^N) \right\},\,  A_V \varphi := -\frac{1}{2}\Delta \varphi + V(x) \varphi$$
For appropriate potentials $V$, approximate controllability of (\ref{Schro}) in $\mathcal{S}$ (possibly in large time) is
a corollary of a general result by Boscain, Caponigro, Chambrion, Mason and Sigalotti
(the original proof of \cite{Chambrion-et-al} is generalized in \cite{Boscain_CCS};
inequality (\ref{T:LB}) below is proved in \cite[Proposition 4.6]{Chambrion-et-al};
see also \cite{Boscain_open_pb} for a survey of results in this area).

\begin{thm} \label{thm:Boscain}
We assume that
\begin{itemize}
\item there exists a Hilbert basis $(\phi_k)_{k \in \mathbb{N}}$ of $L^2(\mathbb{R}^N,\mathbb{C})$ made of eigenvectors of $A_V$:
$A_V \phi_k = \lambda_k \phi_k$ and $x \phi_k \in L^2(\mathbb{R}^N)$, $\forall k \in \mathbb{N}$,
\item $\int_{\mathbb{R}^N} x \phi_j(x) \phi_k(x) dx = 0$ for every $j, k \in  \mathbb{N}$ such that $\lambda_j=\lambda_k$ and $j \neq k$,
\item 
for every $j,k \in \mathbb{N}$, there exists a finite number of integers $p_1,...,p_r \in \mathbb{N}$ such that
\end{itemize}
$$\begin{array}{c}
p_1=j, \quad p_{r}=k, \quad \int_{\mathbb{R}^N} x \phi_{p_l}(x)\phi_{p_{l+1}}(x) dx  \neq 0, \forall l=1,...,r-1\,,
\\
|\lambda_L-\lambda_M| \neq |\lambda_{p_l}-\lambda_{p_{l+1}}|, \forall  1\leqslant l \leqslant r-1, L,M \in \mathbb{N} \text{ with } \{L,M\} \neq \{p_l,p_{l+1}\}.
\end{array}$$

Then, for every $\epsilon>0$ and $\psi_0, \psi_f \in \mathcal{S}$,
there exist a time $T>0$ and a piecewise constant function $u:[0,T] \rightarrow \mathbb{R}$
such that the solution of (\ref{Schro}) satisfies
\begin{equation}\label{presdubut}
  \| \psi(T) - \psi_f \|_{L^2(\mathbb{R}^N)}<\epsilon\,.
\end{equation}
Moreover, for every $\delta >0$,  the existence of a piecewise constant function $u:[0,T] \rightarrow  (-\delta,\delta)$
 such that the solution of (\ref{Schro}) satisfies \eqref{presdubut} implies that
\begin{equation} \label{T:LB}
T \geqslant \frac{1}{\delta} \underset{k \in \mathbb{N}}{\sup}
\frac{|\left|\langle \phi_k,\psi_0\rangle|-|\langle\phi_k,\psi_f\rangle \right||-\epsilon}{\|B\phi_k\|}.
\end{equation}
\end{thm}

To prove this statement, the authors use finite dimensional techniques applied to the Galerkin approximations of equation (\ref{Schro}).
They also prove an estimate on the $L^1$-norm of the control \cite[Proposition 2.8]{Boscain_CCS} and
approximate controllability in the sense of density matrices \cite[Theorem 2.11]{Boscain_CCS}.


In Theorem \ref{thm:Boscain}, the  time $T$ is not known a priori and may be large.
Note that the lower bound on the control time in (\ref{T:LB}) goes to zero when $\delta \rightarrow + \infty$.
Thus, approximate controllability \emph{in arbitrarily small time} (allowing potentially large controls) is an open problem.
The goal of this article is to prove that, for potentials $V$ satisfying (\ref{Ass_V}), 
approximate controllability does not hold in arbitrarily small time, even with large controls,
as stated in the following theorem.

\begin{thm} \label{thm:Main}
Consider $V$ satisfying assumption (\ref{Ass_V}).
Let $b>0$, $x_0, \dot{x}_0 \in \mathbb{R}^N$ and $\psi_0 \in \mathcal{S}$ be defined by
\begin{equation} \label{def:psi0}
\psi_0(x):=\frac{b^{N/4}}{C_N} e^{-\frac{b}{2} \|x-x_0\|^2 + i \langle \dot{x}_0 , x-x_0 \rangle}
\end{equation}
where
$$C_N:=\left(\int_{\mathbb{R}^N} e^{-\|y\|^2} dy \right)^{1/2}\,.$$
Let $\psi_f \in \mathcal{S}$ a state that does not have a Gaussian profile in the sense that
$$|\psi_f(.)| \neq \frac{\det(S)^{1/4}}{C_N} e^{-\frac{1}{2}\| \sqrt{S}(.-\gamma)\|^2}\,, \quad
\forall \gamma \in \mathbb{R}^N, S \in \mathcal{M}_N(\mathbb{R}) \text{ symmetric positive}.$$
Then there exist $T^{**}=T^{**}(\|V''\|_\infty,\|V^{(3)}\|_\infty, b,\psi_f)>0$ and $\delta=\delta(\|V''\|_\infty, b,\psi_f)>0$ such that,
for every $E \in C^0_{pw}([0,T],\mathbb{R}^N)$ (piecewise continuous functions $[0,T]\to \R ^N$), the solution
$\psi $ of (\ref{Schro}) satisfies
$$\| \psi(t)-\psi_f \|_{L^2(\mathbb{R}^N)} > \delta\,, \quad \forall t \in [0,T^{**}]\,.$$
\end{thm}

In particular, if $V$ satisfies (\ref{Ass_V}) and the assumptions of Theorem \ref{thm:Boscain} 
(which hold generically, this fact may be proved as in \cite{Mason_Sigalotti}), then  system (\ref{Schro})
is approximately controllable in $\mathcal{S}$ in large time but not in small time $T<T^{**}$.
In this sense, Theorem \ref{thm:Boscain} is optimal with respect to the time of control. 
A characterization of the minimal time required for $\epsilon$-approximate controllability is an open problem.

\subsection{Bibliographical comments}

\subsubsection{Small-time control and minimal time for ODEs and PDEs}

In \cite{DAlessandro}, D'Alessandro considers (generalizations of) Schr\"{o}dinger ODEs and the controllability of their resolvent.
He characterizes the set of states reachable in arbitrary time from the identity of the group.
In particular, this set may not be the whole compact matrix Lie group even if the system is controllable.
Then, in \cite{Agrachev_Chambrion}, Agrachev and Chambrion prove an estimation of the minimal time for the global approximate controllability.
Such an estimate for PDE (\ref{Schro}) is a widely open problem.

In \cite{Boussaid_Caponigro_Chambrion}, Boussa\"{\i}d, Caponigro and Chambrion present an example of bilinear conservative system in infinite dimension
for which approximate controllability holds in arbitrary small times. This situation is in contrast with the finite dimensional case discussed above
and with Theorem \ref{thm:Main}.

\subsubsection{Minimal time for local exact controllability with small controls}

A different notion of `minimal time' is investigated in \cite{KB-Morancey}.
This article focuses on exact controllability and small controls to realize small motions 
whereas the present article investigates approximate controllability, large controls $E$ and large motions.
Generalizing \cite{JMC-CRAS-Tmin}, the authors of \cite{KB-Morancey} describe a general scenario 
 for local exact controllability (with small controls)
to hold in large time, but not in small time. This positive minimal time is 
related to the loss of directions of the linearized system and the behaviour of the second order term in the power series expansion of the solution.

\subsection{Notations}

Denote by $\mathcal{M}_N(\mathbb{K})$ the set of $N \times N$ matrices with coefficients in $\mathbb{K}=\mathbb{R}$ or $\mathbb{C}$ and $I_N$ its identity element; 
$\text{Tr}(M)$ the trace of a matrix $M \in \mathcal{M}_N(\mathbb{C})$; 
$\mathcal{S}_N(\mathbb{R})$ (resp. $\mathcal{S}_N^+(\mathbb{R})$) the set of symmetric matrices (resp. positive symmetric matrices) in $\mathcal{M}_N(\mathbb{R})$; 
$A \leqslant B$ when $A,B \in \mathcal{S}_N(\mathbb{R})$ and $B-A \in \mathcal{S}_N^+(\mathbb{R})$;
$\|.\|$ the Euclidean norm on $\mathbb{R}^N$ and the associated operator norm on $\mathcal{M}_N(\mathbb{R})$;
$\dot{x}(t):=\frac{dx}{dt}(t)$, $\Ddot{x}(t):=\frac{d^2 x}{dt^2}(t)$, for a function $x$ of the variable $t$;
$C^0_{pw}([0,T],\mathbb{R}^N)$ the piecewise continuous functions $[0,T] \to \R ^N$;
$\langle x, y\rangle :=\sum_{i=1}^N x_iy_i$,
 for every $x=(x_1,x_2,\ldots, x_n)\in \mathbb{C}^N$,
$y=(y_1,y_2,\ldots, y_n)\in \mathbb{C}^N$;
and $\mathcal{S}$ the unit sphere in $L^2(\mathbb{R}^N,\mathbb{C})$.

\section{Proof of Theorem~\ref{thm:Main}}

Let $V$ satisfying (\ref{Ass_V}), $b>0$, $x_0, \dot{x}_0 \in \mathbb{R}^N$ and $\psi_0$ defined by (\ref{def:psi0}).
\\

Our strategy to prove Theorem \ref{thm:Main}, outlined in \cite{Holger_JMP}, is semi-classical: 
it relies on Gaussian approximate solutions that are localized around classical trajectories.
They are called `trajectory-coherent states'  (TCS) and were originally introduced by Bagrov at al. \cite{BBT83,BK92,BK94,BBT96} 
(for recent and comprehensive mathematical treatments, see  \cite{robert2004,rober2012}).
They generalize the well-known explicit solutions for the harmonic oscillator potential $V(x)=x^2$
(see e.g. \cite{HTC,CLT,HarmOsc}) and may also be viewed as generalized WKB states.

The approximate solutions $\widetilde{\psi} = \widetilde{\psi} (t,x)$ (defined in eq. \eqref{def:psitilde} below)
depend on functions $x_c:\R \to \R^N$ and $Q:\R \to \mathcal{M} _N(\C )$, which satisfy the ODEs \eqref{eq:Newton} below.
The vector function $x_c(t)$ is the classical (controlled) trajectory satisfying Newton's equation of motion \eqref{eq:Newton},
which includes the control field $E(t)$.

The remainder of this section is organized as follows. In Section \ref{subsec:Prl},
we prove a preliminary result for the solutions $Q(t)$ of \eqref{eq:Newton}. 
In Section \ref{subsec:approx_sol}, we introduce the explicit approximate solution $\widetilde{\psi}$ and prove that the error
   $\|\psi-\widetilde{\psi}\|_{L^\infty((0,T),L^2(\mathbb{R}^N))}$ can be bounded uniformly with respect to
   $E \in C^0_{pw}(\mathbb{R},\mathbb{R}^N)$.
Finally, Section \ref{subsec:proof} contains the proof of Theorem \ref{thm:Main}.

\subsection{The ODE for $Q(t)$}
\label{subsec:Prl}

For  $E \in C^0_{pw}(\mathbb{R},\mathbb{R}^N)$, we introduce
the maximal solutions $x_c \in C^1 \cap C^2_{pw}(\mathbb{R},\mathbb{R}^N)$ 
and $Q \in C^1((T^-,T^+),\mathcal{M}_N(\mathbb{C}))$ of
\begin{equation}\label{eq:Newton}
\begin{array}{ll}
\left\lbrace \begin{array}{l}
\frac{d^2 x_c}{dt^2}(t)+ \nabla V[x_c(t)]=E(t)\,, \\
x_c(0)=x_0\,, \\
\frac{dx_c}{dt}(0)=\dot{x}_0\,,
\end{array}\right.
&
\left\lbrace \begin{array}{l}
\frac{dQ}{dt}(t) + Q(t)^2 + V''[x_c(t)]=0\,, \\
Q(0)=ib I_N\,,
\end{array}\right.
\end{array}
\end{equation}
where $\nabla V$ and $V''$ denote the gradient and Hessian matrix of $V$, respectively.
Note that $x_c$ is defined for every $t\in \mathbb{R}$ because $\nabla V$ is globally Lipschitz by assumption (\ref{Ass_V});
and the complex coefficient matrix $Q(t)$ is symmetric for every $t \in (T^-,T^+)$,
so $Q_2(t):=\Im[Q(t)]$ is symmetric as well. 
A priori, the maximal interval $(T^-,T^+)$ may depend on $E$.

\begin{Prop} \label{Prop:EDO}
There exists $T^*=T^*(b,\|V''\|_\infty)>0$ such that, for every $E \in C^0_{pw}(\mathbb{R},\mathbb{R}^N)$,
$Q(t)$ is defined for every $t \in [0,T^{*}]$ (i.e. $T^+>T^*$) and
\begin{equation}\label{ineqQ2}
\frac{b}{2} I_N \leqslant Q_2(t) \leqslant \frac{3b}{2} I_N\,,\text{ for every } t \in [0,T^*]\,.
\end{equation}
\end{Prop}
\noindent \textbf{Proof of Proposition \ref{Prop:EDO}:} Let $T^*=T^*(b,\|V''\|_\infty)>0$ be such that
\begin{equation} \label{def:T*}
t[1+be^{4t}+\|V''\|_\infty]<1\, \quad \text{ and } \quad 2 t e^{2t} \leqslant \frac{1}{2}\,, \quad  \forall t \in [0,T^*]\,.
\end{equation}

\textbf{Step 1: Equations satisfied by $Q_1(t):=\Re[Q(t)]$ and $Q_2(t):=\Im[Q(t)]$.}
Since $V''$ is real, (\ref{eq:Newton}) implies, on $(T^-,T^+)$
\begin{equation} \label{eq:Q2}
\begin{array}{ll}
\left\lbrace \begin{array}{l}
\frac{dQ_1}{dt} + Q_1^2 - Q_2^2 + V''[x_c]=0\,,  \\
Q_1(0)=0\,,                                        
\end{array}\right.
&
\left\lbrace \begin{array}{ll}
\frac{dQ_2}{dt} + Q_1 Q_2 + Q_2 Q_1 =0\,, \\
Q_2(0)=b I_N\,.                                      
\end{array}\right.
\end{array}
\end{equation}
By Gronwall lemma and \eqref{eq:Q2},
\begin{equation} \label{Q2:explicit}
\|Q_2(t)\| \leqslant b e^{2 \int_0^t \|Q_1(s)\| ds}\,, \quad \forall t \in [0,T^+)\,.
\end{equation}

\textbf{Step 2: We prove that $T^+ > T^*$ for every $E \in C^0_{pw}(\mathbb{R},\mathbb{R}^N)$.}
Working by contradiction, we assume the existence of $E \in C^0_{pw}(\mathbb{R},\mathbb{R}^N)$ such that $T^+ < T^*$.
In particular $T^+$ is finite, thus $Q(t)$ explodes as $t \rightarrow T^+$.
We then deduce from (\ref{Q2:explicit}) that $Q_1(t)$ explodes as $t \rightarrow T^+$.
Thus
\begin{equation}\label{deft*}
  t^*:=\sup\{ t \in [0,T^+) ; \|Q_1(s)\| \leqslant 1\,, \forall s \in [0,t] \}
\end{equation}
belongs to $(0,T^+)$ and $\|Q_1(t^*)\|=1$. Now we have
$$\begin{array}{ll}
1 & = \|Q_1(t^*)\|  = \Big\| \int_0^{t^*} \Big(- Q_1(s)^2 + Q_2(s)^2 + V''[x_c(s)]\Big) ds  \Big\| \quad \text{ by } (\ref{eq:Q2})\\
  & \leqslant \int_0^{t^*} \Big( \|Q_1(s)\|^2 + b^2 e^{4 \int_0^t \|Q_1(s)\|ds}  + \|V''[x_c(s)]\| \Big) ds  \quad \text{ by } (\ref{Q2:explicit})\\
  & \leqslant t^* [ 1 + b^2 e^{4t^*}+\|V''\|_\infty ] < 1 \quad \text{ by } (\ref{deft*}) \text{ and } (\ref{def:T*}),
\end{array}$$
which is a contradiction. 
Therefore $T^+ > T^*$ for every $E \in C^0_{pw}(\mathbb{R},\mathbb{R}^N)$.

\textbf{Step 3: Conclusion.} The same argument proves that $t^* > T^*$, i.e. $\|Q_1(t)\| \leqslant 1$ for every $t \in [0,T^*]$.
Thus, by \eqref{eq:Q2}, (\ref{Q2:explicit}) and (\ref{def:T*}) we have
$$\| Q_2(t)- b I_N \| = \left\|  \int_0^t [ Q_2 Q_1 + Q_1 Q_2 ](s) ds \right\|  \leqslant 2 t b e^{2t} < \frac{b}{2}\,, \forall t \in [0,T^*]\,.$$

\subsection{Approximate solution}
\label{subsec:approx_sol}

We introduce the `classical action' $S:(t,x) \in \mathbb{R} \times \mathbb{R}^N \rightarrow \mathbb{R}$
\begin{equation}\label{defS}
  S(t,x):= \int_0^t \left( \frac{1}{2} \|\dot{x}_c(s)\|^2 - V[x_c(s)] \right) ds + \langle \dot{x}_c(t) , x-x_c(t)\rangle
\end{equation}
and the approximate solution
\begin{equation} \label{def:psitilde}
\begin{array}{c}
\widetilde{\psi}(t,x):= \frac{b^{N/4}}{C_N} \exp\Big[ \Phi(t,x) \Big] \text{ where } \\
\begin{array}{ll}
\Phi(t,x):= & i \Big( S(t,x)+\frac{1}{2}\langle Q(t) [x-x_c(t)],x-x_c(t)\rangle \Big) +\\
            &  \int_0^t \left(i \langle x_c(s),E(s)\rangle - \frac{\text{Tr}[Q(s)]}{2} \right) ds  \Big]\,.
\end{array}
\end{array}
\end{equation}
Note that equation (\ref{eq:Q2}) ensures that $\widetilde{\psi}(t) \in \mathcal{S}$. Indeed, for every $t\in (0,T^*)$,
$$\begin{array}{ll}
\frac{d}{dt} \det[Q_2(t)] & = \det[Q_2(t)] \text{Tr}\left[ Q_2(t)^{-1} \dot{Q}_2(t) \right] \\
                          & = \det[Q_2(t)] \text{Tr}\left[ Q_2(t)^{-1} \Big( - Q_1(t)Q_2(t)-Q_2(t)Q_1(t) \Big) \right] \\
                          & = - 2 \text{Tr}[Q_1(t)] \det[Q_2(t)]\,,
\end{array}$$
which, together with \eqref{eq:Q2}, implies that
\begin{equation} \label{det(Q2):explicit}
\det[Q_2(t)] = b^N e^{-2\int_0^t \text{Tr}[Q_1(s)]} ds\,.
\end{equation}
Therefore,
$$\begin{array}{ll}
\|\widetilde{\psi}(t)\|_{L^2(\mathbb{R}^N)}
& = \displaystyle \frac{b^{N/4}}{C_N} e^{- \frac{1}{2} \int_0^t \text{Tr}[Q_1(s)]ds} \left( \int_{\mathbb{R}^N} e^{-\langle Q_2(t)[x-x_c(t)],x-x_c(t)\rangle } dx \right)^{1/2}\\
& = \displaystyle \frac{b^{N/4}}{C_N} e^{- \frac{1}{2} \int_0^t \text{Tr}[Q_1(s)]ds} \frac{C_N}{\det[Q_2(t)]^{1/4}}   =1\,.
\end{array}$$

\begin{Prop} \label{Prop:error}
There exists a constant $C_{*}>0$ such that,
for every $V$ satisfying (\ref{Ass_V}), $b>0$, $x_0, \dot{x}_0 \in \mathbb{R}^N$ and $E \in C^0_{pw}(\mathbb{R},\mathbb{R}^N)$,
the solution $\psi$ of (\ref{Schro}) with $\psi_0$ defined by (\ref{def:psi0})
and the function $\widetilde{\psi}$ defined by (\ref{def:psitilde}) satisfy
$$\| (\psi-\widetilde{\psi})(t) \|_{L^2(\mathbb{R}^N)} \leqslant
C_{*} \| V^{(3)}\|_{\infty} \int_0^t \|Q_2(s)^{-1}\|^{3/2} ds \,, \quad \forall t \in [0,T^*]\,,$$
where $T^*$ is defined in Proposition \ref{Prop:EDO}.
\end{Prop}

\noindent \textbf{Proof of Proposition \ref{Prop:error}:}

\textbf{Step 1: Equation satisfied by $\widetilde{\psi}$.}
On the one hand, we have (note that $S$ is linear in the variable $x$)
$$\nabla \widetilde{\psi} = i \Big( \nabla S + Q [x-x_c] \Big)  \widetilde{\psi}\,,$$
$$\begin{array}{ll}
\Delta \widetilde{\psi}
& = \left[  -  \langle \nabla S+ Q [x-x_c],\nabla S+ Q [x-x_c] \rangle
 + i \text{Tr}(Q) \right]  \widetilde{\psi}   \\
& = \left[  - \langle \dot{x}_c + Q [x-x_c],\dot{x}_c + Q [x-x_c]\rangle  + i \text{Tr}(Q) \right]  \widetilde{\psi}\,. \\
\end{array}$$
Thus
\begin{equation} \label{dx2}
\frac{1}{2} \Delta \widetilde{\psi} =
\left[
-\frac{1}{2} \|\dot{x}_c\|^2 - \langle \dot{x}_c , Q [x-x_c]\rangle - \frac{1}{2}
\langle Q^2 [x-x_c],[x-x_c]\rangle  + \frac{i}{2} \text{Tr}(Q)
\right] \widetilde{\psi}\,.
\end{equation}
On the other hand, we have
\begin{equation} \label{idt}
\begin{array}{lll}
i \partial_t \widetilde{\psi}
& = &
\left[ -\Big( \partial_t S +\frac{1}{2}\langle \dot{Q}[x-x_c],x-x_c\rangle - \langle Q \dot{x}_c, x-x_c \rangle \Big) - \frac{i}{2}\text{Tr}(Q) - \langle x_c, E\rangle \right] \widetilde{\psi}
\\  & = &
\Big[
- \frac{1}{2} \|\dot{x}_c\|^2 + V[x_c] - \langle \Ddot{x}_c , x-x_c\rangle + \|\dot{x}_c\|^2
- \frac{1}{2}\langle \dot{Q} [x-x_c],x-x_c\rangle
\\ & & + \langle Q \dot{x}_c , x-x_c \rangle -  \frac{i}{2} \text{Tr}(Q) - \langle x_c , E \rangle  \Big] \widetilde{\psi}
\\ & = &
\Big[
\frac{1}{2} \|\dot{x}_c\|^2 + V[x_c] + \langle \nabla V(x_c)-E , x-x_c \rangle
\\ & &
+ \frac{1}{2} \langle [Q^2+V''(x_c)][x-x_c],x-x_c\rangle  + \langle Q \dot{x}_c , x-x_c \rangle
\\ & & -  \frac{i}{2} \text{Tr}(Q) - \langle x_c , E \rangle
\Big] \widetilde{\psi}\,.
\end{array}
\end{equation}
Combining (\ref{dx2}) and (\ref{idt}) gives
\begin{equation} \label{eq:psi_tilde}
i \partial_t \widetilde{\psi}(t,x) + \frac{1}{2} \Delta  \widetilde{\psi}(t,x) - V(x) \widetilde{\psi}(t,x) + \langle E(t) , x \rangle \widetilde{\psi}(t,x)
= r(t,x)
\end{equation}
where
\begin{equation} \label{def:r}
r(t,x):=- \Big( V(x) - V(x_c) - \langle \nabla V(x_c) , x-x_c \rangle - \frac{1}{2} \langle V''(x_c)[x-x_c],x-x_c\rangle \Big) \widetilde{\psi}(t,x)\,.
\end{equation}

\textbf{Step 2: Conclusion.}
Using \eqref{def:r}, (\ref{def:psitilde}) and (\ref{det(Q2):explicit}) and Taylor's formula, we get
\begin{eqnarray*}
\|r(t)\|_{L^2(\mathbb{R}^N)}^2
& \leqslant & \int_{\mathbb{R}^N} \Big| \frac{\|V^{(3)}\|_\infty}{3!} \|x-x_c\|^3 \Big|^2 \frac{b^{N/2}}{C_N^2} e^{-\langle Q_2(t)[x-x_c],x-x_c\rangle -\int_0^t \text{Tr}(Q_1)} dx \\
& \leqslant &\frac{\| V^{(3)}\|_{\infty}^2}{(3!)^2 C_N^2}  \int_{\mathbb{R}^N} \|x-x_c\|^6  e^{-\langle Q_2(t)[x-x_c],x-x_c\rangle} \sqrt{\det[Q_2(t)]} dx \\
& \leqslant & C_{*}^2 \| V^{(3)}\|_{\infty}^2 \|Q_2(t)^{-1/2}\|^6\,,
\end{eqnarray*}
where
$$C_{*}:=\frac{1}{3! C_N} \left(\int_{\mathbb{R}^N} \|y\|^6 e^{-\|y\|^2} dy\right)^{1/2}\,.$$
Let $U(t,s)$ be the evolution operator for equation (\ref{Schro}) (see Proposition \ref{Prop:WP}). Then,
$$(\psi-\widetilde{\psi})(t)=\int_0^t U(t,s) r(s) ds \text{ in } L^2(\mathbb{R}^N)\,, \forall t \in (0,T^*)\,,$$
and $U(t,s)$ is an isometry of $L^2(\mathbb{R}^N)$ for every $t \geqslant s \geqslant 0$, thus
\begin{equation}\label{estpsi-tildepsi}
\|(\psi-\widetilde{\psi})(t)\|_{L^2(\mathbb{R}^N)}
\leqslant \int_0^t \|r(s)\|_{L^2(\mathbb{R}^N)} ds
\leqslant \int_0^t  C_{*} \| V^{(3)}\|_{\infty} \|Q_2(s)^{-1}\|^{3/2} ds.
\end{equation}
\hfill $\Box $

\begin{rk}
The trajectory-coherent states $\widetilde{\psi}$ are only approximate solutions to the Schr\"odinger equation;
however they are exact solutions for quadratic potentials $V$ (see (\ref{eq:psi_tilde}) and (\ref{def:r})),
which is the key point of reference \cite{HarmOsc}.
\end{rk}

\subsection{Proof of the main result}
\label{subsec:proof}

Let $T^*=T^*(b,\|V''\|_\infty)>0$ be as in Proposition \ref{Prop:EDO}.
The key point of the proof is the fact that, for every $E \in C^0_{pw}(\mathbb{R},\mathbb{R})$,
the approximate solution $\widetilde{\psi}$ has a Gaussian profile. Indeed,
from  \eqref{defS}, \eqref{def:psitilde} and \eqref{det(Q2):explicit}, one has
$$|\widetilde{\psi}(t,x)|^2 =  \frac{\det[\sqrt{Q_2(t)}]}{C_N^2} e^{-\left\|\sqrt{Q_2(t)}[x-x_c(t)]\right\|^2  }\,,$$
where $Q_2(t)$ is a real symmetric matrix satisfying \eqref{ineqQ2}.

\textbf{Step 1: We prove that the set}
$$\begin{array}{ll}
\mathcal{V}:= & \Big\{ \phi \in \mathcal{S} ;
\exists q \in \mathcal{S}_N^+(\mathbb{R}) \text{ with } \sqrt{\frac{b}{2}} I_N \leqslant q \leqslant \sqrt{\frac{3b}{2}} I_N \text{ and }
\\ &\quad  \alpha \in \mathbb{R}^N \text{ such that } |\phi(x)|^2 = \frac{\det(q)}{C_N^2} e^{- \|q(x-\alpha)\|^2} \text{ a.e. }\Big\}
\end{array}$$
\textbf{is a strict closed subset of $\mathcal{S}$.}

Clearly, $\mathcal{V}$ is a strict subset of $\mathcal{S}$.
Let $(\phi_n)_{n \in \mathbb{N}}$ be a sequence of $\mathcal{V}$ that converges in $L^2(\mathbb{R}^N,\mathbb{C})$ to $\phi_\infty \in \mathcal{S}$.
For every $n \in \mathbb{N}$, we denote by $\alpha_n \in \mathbb{R}^N$ and $q_n \in \mathcal{S}_N^+(\mathbb{R})$
the corresponding parameters;
$q_n$ satisfies
\begin{equation} \label{ineq_qn}
\sqrt{\frac{b}{2}} I_N \leqslant q_n \leqslant \sqrt{\frac{3b}{2}} I_N.
\end{equation}
By extracting a subsequence if necessary,  we may assume w.l.o.g.
that $\phi_n(x) \rightarrow \phi_\infty(x)$ for almost every $x \in \mathbb{R}^N$ (Lebesgue).

\textbf{Step 1.1: Up to a possible extraction of a subsequence, we may assume that $q_n \rightarrow q_\infty$} where
$q_\infty \in \mathcal{S}_N^+(\mathbb{R})$ and $\sqrt{\frac{b}{2}} I_N \leqslant q_\infty \leqslant \sqrt{\frac{3b}{2}} I_N$.
This may be seen by diagonalizing and appealing to the compactness of $O_N(\mathbb{R}) \times [\sqrt{b/2},\sqrt{3b/2}]^N$.

\textbf{Step 1.2: Up to a possible extraction of a subsequence, we may assume that $\alpha_n \rightarrow \alpha_\infty\in \mathbb{R}^N$.}
Working by contradiction, we assume that $(\alpha_n)_{n \in \mathbb{N}}$ is not bounded.
We may assume w.l.o.g. that $\|\alpha_n\| \rightarrow + \infty$ when $n \rightarrow \infty$.
Then, by (\ref{ineq_qn}),
$$|\phi_n(x)|^2
= \frac{\det(q_n)}{C_N^2} e^{- \|q_n(x-\alpha_n)\|^2}
\leqslant \frac{(3b/2)^{N/2}}{C_N^2} e^{-\frac{b}{2}\|x-\alpha_n\|^2} \rightarrow 0 \quad \text{ a.e. } x \in \mathbb{R}^N\,.$$
Thus $\phi_\infty=0$ (uniqueness of the a.e. limit), which is impossible because $\phi_\infty \in \mathcal{S}$.

\textbf{Step 1.3: Conclusion.}
The uniqueness of the a.e. limit gives
$$|\phi_\infty(x)|^2 =\frac{\det(q_\infty)}{C_N^2} e^{-\| q_\infty(x-\alpha_\infty)\|^2} \quad \text{ a.e. } x \in \mathbb{R}^N$$
thus $\phi_\infty \in \mathcal{V}$. This concludes Step 1.
\\

\textbf{Step 2:} Let $\psi_f \in \mathcal{S} \setminus \mathcal{V}$
(which holds, in particular, when $\psi_f$ does not have a Gaussian profile).
Then
$\delta_0:=\text{distance}_{L^2(\mathbb{R}^N,\mathbb{C})}( \psi_f ; \mathcal{V} )>0$.
Let
$$T^{**}=T^{**}(\psi_f,b,V):=\min\left\{ T^* ; \frac{\delta_0 (b/2)^{3/2}}{2 C_{*} \|V^{(3)}\|_\infty} \right\}\,.$$
Then, using \eqref{ineqQ2} and \eqref{estpsi-tildepsi}, we get that, for every $t \in [0,T^{**}]$ and $E \in C^0_{pw}(\mathbb{R},\mathbb{R})$,
the solution $\psi$ of (\ref{Schro}) satisfies
$$\begin{array}{ll}
\|\psi_f - \psi(t)\|_{L^2(\mathbb{R})} & \geqslant \Big| \|\psi_f - \widetilde{\psi}(t)\|_{L^2(\mathbb{R})} - \| (\widetilde{\psi}-\psi)(t) \|_{L^2(\mathbb{R})} \Big| \\
                                       &\displaystyle  \geqslant \delta_0 - C_{*} \|V^{(3)}\|_\infty \frac{t}{(b/2)^{3/2}} \geqslant \frac{\delta_0}{2}\,. \hfill \Box
\end{array}$$

\textbf{Acknowledgements:} The authors were partially supported by the ``Agence Nationale de la Recherche'' 
(ANR) Projet Blanc EMAQS number ANR-2011-BS01-017-01 and by the ERC advanced grant 266907
(CPDENL) of the 7th Research Framework Programme (FP7).

\bibliography{biblio4}
\bibliographystyle{plain}

\end{document}